\documentclass[3p]{elsarticle}

\usepackage{graphicx,amsmath,amssymb,txfonts}
\usepackage{amsmath,amssymb}
\usepackage{algorithm}
\usepackage{algorithmic}
\usepackage{subfigure,color}
\usepackage[colorlinks]{hyperref}
\usepackage{hypernat}
\usepackage{multirow}
\usepackage{booktabs}
\usepackage{threeparttable}
\usepackage{epstopdf}

\graphicspath{{Figs/}}

\makeatletter

\newcommand{\Rmnum}[1]{\expandafter\@slowromancap\romannumeral #1@}
\makeatother

\newtheorem{remark}{\textbf{Remark}}
\newtheorem{theorem}{Theorem}
\newtheorem{example}{Example}
\begin{document}
\begin{frontmatter}
\title{Fourth-order compact finite difference schemes for solving biharmonic equations with Dirichlet boundary conditions}
 \author[csu]{Kejia Pan}
  \author[csu]{Jin Li\corref{cor1} }
  \author[ncsu]{Zhilin Li}
   \author[csu]{Kang Fu}
 \address[csu]{School of Mathematics and Statistics, Central South University, Changsha 410083, China}
 \address[ncsu]{CRSC \& Department of Mathematics, North Carolina State University, Raleigh, NC 27695-8205, USA}
  \cortext[cor1]{Corresponding authors. \\ kejiapan$@$csu.edu.cn(K.J. Pan) li\_jin$@$csu.edu.cn(J. Li) zhilin$@$math.ncsu.edu(Z.L. Li) fukang$@$csu.edu.cn(K. Fu) }

\begin{abstract}
In this study, we propose a genuine fourth-order compact finite difference scheme for solving biharmonic equations with Dirichlet boundary conditions in both two and three dimensions.
In the 2D case, we build upon the high-order compact (HOC) schemes for flux-type boundary conditions originally developed by Zhilin Li and Kejia Pan [SIAM J. Sci. Comput., 45 (2023), pp. A646-A674] to construct a high order compact discretization for coupled boundary conditions.
When considering the 3D case, we modify carefully designed undetermined coefficient methods of Li and Pan to derive the finite difference approximations of coupled boundary conditions. 
The resultant FD discretization maintains the global fourth order convergence and compactness.
Unlike the very popular Stephenson method, the number of unknows do not increase with dimensions.
Besides, it is noteworthy that the condition number of the coefficient matrix increases at a rate of $O(h^{-2})$ in both 2D and 3D. 
We also validate the performance of the proposed genuine HOC methods through nontrivial examples.

\end{abstract}

\begin{keyword}
Biharmonic problems\sep Fourth order compact \sep Finite difference method \sep Coupled Poisson system \sep Stokes flow
\MSC 65N06 \sep 65N55
\end{keyword}

\end{frontmatter}

\section{Introduction}
The biharmonic equation, a fourth order partial differential equation(PDE), plays a pivotal role in various branches of applied mathematics and engineering, particularly in the fields of elasticity theory, fluid mechanics, and plate bending problems. This equation arises naturally when modeling physical phenomena that involve the bending of elastic plates or the flow of viscous fluids. Given its wide range of applications, the accurate and efficient solution of the biharmonic equation is of great interest in both theoretical and computational studies.

In many practical scenarios, the biharmonic equation is subject to boundary conditions that significantly influence the behavior of its solutions. Among these, Dirichlet boundary conditions are particularly common, where the values of the function and its derivatives are specified on the boundary of the domain. The accurate enforcement of these boundary conditions is crucial for obtaining physically meaningful solutions. However, the complexity of the biharmonic equation, compounded by the stringent nature of Dirichlet boundary conditions, presents substantial challenges in deriving accurate and efficient numerical solutions.
 
 For biharmonic problems, a fourth-order PDE, there are two categories of numerical approaches that are widely applied. 
 The first one is called a direct or uncoupled approach in which considers how to approximate the fourth order differential operator, found in \cite{FD2d13,FD3d27,SFD1984,SFD1998,SFD2008}. 
 Another one is named the coupled, splitting, or mixed approach, which introduces the Laplace of primary unknowns as intermediate unknowns such that a fourth-order problem can be reduced to two two-order PDEs, referring \cite{CFD1968,CFD1970,CFD1971_1,CFD1971_2,CFD1997}.
 
 The direct approach combined with finite difference (FD) discretization also fall broadly into two categories. 
 Here, \cite{FD2d13} and \cite{FD3d27} directly approximate the biharmonic operator with second-order accuracy but use 13-points in 2D and 27-points in 3D, which destroy the compactness of the FD method and need extra treatments near the boundary.
 In \cite{SFD1984}-\cite{SFD1984}, a Stephenson method (also called the combined FD method) and its modifications are proposed to discretize biharmonic problems with second- or fourth-order accuracy.
  The method uses the gradients of the primary unknown, even second-order derivatives, for high-order accuracy as intermediate unknowns, which results in the number of unknowns increasing with the dimension. 
  Such as, considering a second-order accuracy, on a $N\times N$ grid, the degree of freedom of Stephenson schemes is $3 N^2$, while on a $N\times N\times N$ grid in 3D, the degree of freedom increases to $4N^3$.
  
 The coupled approach aims to reduce the biharmonic problem into two Poisson problems. 
 After the treatment, the coupled system is also equal to the Stokes problem in stream-vorticity form. 
 The Stephenson schemes \cite{NA2012} also can be used for such problems, which can be seen as a combined compact FD scheme of Poisson problems. But it leads to the number of unknonws increase with dimensions.
Another straightforward idea is that, through the coupled approach, it is possible finite difference methods or thier high-order formulations for Poisson equations can be applied to the discretization without extra unknowns.
If the Dirichlet boundary conditions of the second are prescribed, the answer is yes. 
But for the Dirichlet boundary condition of the first, the Poisson problems of intermediate unknowns lack suitable boundary conditions, while another is overdetermined.  
There seems to be no research to obtain genuine high-order compact finite difference schemes (both interior and boundary) from what we have seen so far. 
Most of the research chose non-compact or low-order discretization on or near boundaries \cite{IJNMF1995,IJNMF1998,IJNMF2009,NMTMA2019}.
  
%  We can refer \cite{CFD1968}\cite{CFD1970}\cite{CFD1971_1}\cite{CFD1971_2}\cite{CFD1997}.

The primary goal of this paper is to introduce a genuine high-order compact finite difference scheme for biharmonic problems using coupled formulations. Unlike Stephenson schemes, where the number of unknowns increases with the dimension, the proposed schemes maintain a consistent number of unknowns at $2N^d$, where $d$ represents the dimension. Notably, the condition number of the discrete coefficient matrix increases at a rate of $O(h^{-2})$, which is even better than the low-order direct finite difference approaches discussed in \cite{FD2d13} and \cite{FD3d27}.

The structure of this paper is as follows. 
In section \ref{sec:Not}, we provide preparatory material on the coupled approach and finite difference method. 
In section \ref{sec:2d}, we introduce the core of the proposed method, focusing on boundary treatment and the application of results from \cite{li-pan-hoc21} in 2D. 
In section \ref{sec:3d}, we extend the approach by carefully modifying the derivation from \cite{li-pan-hoc21}, employing a method of undetermined coefficients to develop the HOC schemes.
In section \ref{sec:NE}, we present numerical experiments that verify the convergence order and condition number and include an application to Stokes flow.

 \section{Notation} \label{sec:Not}

\subsection{Biharmonic equations}
In the present work, we are mainly concerned with the following biharmonic problem:
\begin{equation}
	\Delta^2 u = f, \;\;\text{in}\; \Omega,
\end{equation}
with general boundary conditions, which include both Dirichlet boundary conditions of the first kind and the second kind
\begin{align}
	u &= g_D, \;\;\text{on}\; \partial \Omega,\label{eq_BC_D} \\
	\partial_n u &= g_N, \;\;\text{on}\; \partial \Omega_1,\label{eq_BC_N1}\\
	 \Delta u & = g_L,\;\;\text{on}\; \partial \Omega_2,\label{eq_BC_N2}
\end{align}
where $\Omega$ is an open, bounded Lipschitz set in $\mathbb{R}^d$, $d=2,3$, with boundary $\partial \Omega = \partial \Omega_1 \cup\partial\Omega_2$, the load $f$ is in $C^\alpha(\Omega)$, the primitive solution $u$ is in $C^{\alpha +4}(\Omega)$ and $\partial_n$ denotes the normal derivatives on $\partial \Omega$. For a fourth order scheme, we need to set $\alpha$ equal to 2 for convergent proof. Then, we briefly introduce the coupled approach, which decomposes the original fourth order problem into two Poisson problems.
We first define an intermediate variable $v=\Delta u$, which is also in $C^{\alpha+2}(\Omega)$ and governed by $u$ on $\partial \Omega$. After, we can rewrite the biharmonic problem by
\begin{align}
	\Delta v  = f,\;\;
	\Delta u = v,
\end{align}
with the same boundary conditions (\ref{eq_BC_D}-\ref{eq_BC_N2}).
It is straightforward that biharmonic equations reduce to two Poisson equations with Dirichlet BCs $(u=g_D,v=g_L)$ when only Dirichlet BCs of the second kind $(\Delta u=g_L)$ are prescribed.
But, when considering boundary conditions (\ref{eq_BC_D}-\ref{eq_BC_N2}), it is challenging to derive appropriate boundary conditions for $v$, especially in constructing a high order accurate scheme. Our main contributions of this paper are to overcome it and design a high order compact finite difference scheme, which also results in a well-conditioned coefficient matrix and avoids nonlinear iteration.

\subsection{Fourth order compact scheme for Poisson equations with a Neumann BC}
Numerous high order compact finite difference schemes have been proposed for elliptic PDEs with Dirichlet BCs. Recently, a genuine high order compact scheme \cite{li-pan-hoc21} for flux BCs has been proposed, which is the cornerstone of the proposed scheme for biharmonic problems.

Here, we recall the case of Poisson equations with Neumann BCs,
\begin{align}
	\Delta u &=f, \;\;\text{in}\; \Omega, \\
	\partial_nu&=g_N, \;\;\text{on}\; \partial \Omega_1,\\
	u&=g_D,\;\;\text{on}\; \partial\Omega_2.
\end{align}
If the square $\Omega=[x_l,x_r]\times[y_l,y_r]$ is covered by a uniform grid 
\begin{equation}
	x_i = x_l+ih,\;\;i=0,\cdots,N_x;\qquad y_j=y_l+jh,\;\;j=0,\cdots,N_y, \label{Grid}
\end{equation}
$h=(x_r-x_l)/N_x=(y_r-y_l)/N_y$ and Neumann boundary condition is prescribed on the left boundary $x=x_l$ that is $\partial\Omega_1=\left\{(x_l,y_j)\vert 0<j<N_y \right\}$.

At interior grid points $(x_i,y_j)$, the classical fourth order compact nine-point finite difference scheme can be written as
\begin{equation}
\begin{aligned}
\frac{1}{6h^2}\left(U_{i-1,j-1}+4U_{i,j-1}+U_{i+1,j-1}+4U_{i-1,j}-20U_{ij}+4U_{i+1,j}+U_{i-1,j+1}\right.\\
\left.+4U_{i,j+1}+U_{i+1,j+1} \right) =\frac{1}{12}\left(f_{i,j-1}+f_{i-1,j}+8f_{ij}+f_{i+1,j}+f_{i,j+1} \right) \label{eq-HOC-int},
\end{aligned}	
\end{equation}

where $\{U_{ij}\}$ is the solution value to the Poisson problem at grid points $(x_i,y_j)$ and $f_{ij}=f(x_i,y_j)$.
Then, refer \cite{li-pan-hoc21}, the fourth order compact scheme at a boundary grid point $(x_i,y_j)$ can be written by
\begin{align}
\frac{1}{6h^2}\left(4U_{i,j-1}+2U_{i+1,j-1}-20U_{ij}+8U_{i+1,j}+4U_{i,j+1}+2U_{i+1,j+1} \right) \notag	\\=\frac{1}{12}\left(f_{i,j-1}-f_{i-1,j}+8f_{ij}+3f_{i+1,j}+f_{i,j+1} \right)-\frac{2}{h}(g_N)_{j},\label{eq-HOC-NB}
\end{align}
where $(g_N)_j=g_N(y_j)$ and $i=0;j=1,\cdots,N_y-1$. Note that the term $f_{i-1,j}$ is defined out of the domain $\Omega$. 
And in \cite{li-pan-hoc21}, the value $f_{i-1,j}$ is assumed to be known or obtained by a quadratic extension. The three dimensional cases can also be found in the research \cite{li-pan-hoc21}.

\section{Fourth order compact schemes for Biharmonic problem in two dimension} \label{sec:2d}
After introducing intermediate variables $v=\Delta u$, the solution of biharmonic problems can be transformed into two Poisson problems. That is,
\begin{align}
	\Delta v &= f,  \;\;\text{in}\; \Omega \label{eq-vP}\\
	\Delta u &= v,   \;\;\text{in}\; \Omega \label{eq-uP}\\
	u &= g_D,  \;\;\text{on}\; \partial\Omega \label{eq-uB}\\
	\partial_n u&=g_N,  \;\;\text{on}\; \partial\Omega_1,\label{eq-vB1}\\
	v &= g_L, \;\;\text{on}\; \partial\Omega_2.\label{eq-vB2}
\end{align}
In the following discussion, we assume $\partial \Omega_1={(x_l,y)}$ such that a clearer derivation can be presented. Then, we will respectively detail the discretization of $u$ and $v$ with fourth order accuracy using the same discrete grid \eqref{Grid}.

If we regard $v$ as a source term of $u$'s Poisson equations \eqref{eq-uP}, 
utilizing fourth order compact scheme \eqref{eq-HOC-int} can generate the following discrete Poisson equations at an interior grid point $(x_i,y_j)$,
\begin{equation}
	\begin{aligned}
\frac{1}{6h^2}\left(U_{i-1,j-1}+4U_{i,j-1}+U_{i+1,j-1}+4U_{i-1,j}-20U_{ij}+4U_{i+1,j}+U_{i-1,j+1}\right.\\
\left.+4U_{i,j+1}+U_{i+1,j+1} \right)-\frac{1}{12}\left(V_{i,j-1}+V_{i-1,j}+8V_{ij}+V_{i+1,j}+V_{i,j+1} \right)=0,
\end{aligned}
\end{equation}
where $\{U_{ij}\}$ and $\{V_{ij}\}$ are the approximate value to the primive varibale $u$ and intermediate varibale $v$ at grid points $(x_i,y_j)$.

Furthermore, by Dirichlet boundary conditions \eqref{eq-uB}, the values of $u$ are known on the boundary. Hence, we get the discrete $u$'s Poisson problem 
\begin{equation}
	A_h\mathcal{U}+B\mathcal{V}=F_1, \label{eq-Dis-uPoiss}
\end{equation}
where  the vectors $\mathcal{U}$ and $\mathcal{V}$ represent the set of $\{U_{ij}\}$ and $\{V_{ij}\}$ respectively.
It is easy to check that the matrix $(-A_h)$ is an M-matrix.

For $v$'s Poisson problem, at an interior grid point $(x_i,y_j)$, we also apply the classical fourth order compact finite difference \eqref{eq-HOC-int},
\begin{equation}
	\begin{aligned}
\frac{1}{6h^2}\left(V_{i-1,j-1}+4V_{i,j-1}+V_{i+1,j-1}+4V_{i-1,j}-20V_{ij}+4V_{i+1,j]+V_{i-1,j+1}}\right.\\
\left.+4V_{i,j+1}+V_{i+1,j+1} \right)=\frac{1}{12}\left(f_{i,j-1}+f_{i-1,j}+8f_{ij}+f_{i+1,j}+f_{i,j+1} \right) \label{}.
\end{aligned}
\end{equation}
The challenge is how to give an artificial boundary condition or construct a discrete formulation of $v$'s Poisson problem on the boundary $x=x_l$, which also will not affect the global accuracy of the numerical method.

Based on the difference of the known boundary conditions, we divided all boundary points into two categories to discuss: corner grid points, e.g., $(x_l,y_l)$, and left boundary grid points.

The numerical treatment for the corner points is relatively straightforward since there are two Dirichlet conditions of $u$ in two directions. 
Such as the corner point $(x_l,y_l)$, we can get $u(x_l,y)=g_N(y)$ and $u(x,y_l)=g_N(x)$ holds from the Dirichlet boundary conditions \eqref{eq-uB} of $u$.
 Hence, by analytical differentiating them, the value of $v$ can be computed, $v(x_l,y_l) = \Delta u(x_l,y_l) = \left( \frac{\partial^2u(x,y_l)}{\partial x^2}+\frac{\partial^2 u(x_l,y)}{\partial y^2} \right)\Big\vert_{(x_l,y_l)}$, which is a Dirichlet boundary condtion for $v$.
 \begin{remark}
 If only boundary conditions \eqref{eq-vB1} are prescribed on four sides, in the same way, the values of $v$ at the other three corner points can be obtained as well. 
 In other words, the variables $v$ are prescribed by Dirichlet BCs at four corner points by the derivation.
 \end{remark}

 Then, we discuss the discretization of other boundary grid points. 
 Take $\left\{ (x_l,y_j)\right.$ $\left.\big\vert j=1,\cdots,N_y-1 \right\}$ as an example. 
 We recall the existing conditions.
 By $u$'s Poisson equation, we first know $v$ is governed by $\Delta u$. Then, the Neumann boundary conditions \eqref{eq-vB1} also are not used. Here, we try to apply the fourth order compact FD scheme \eqref{eq-HOC-NB} for the Neumann BC where $v$ is considered as a source term of $u$'s Poisson equations. That is, 
 \begin{align}
\frac{1}{6h^2}\left(4U_{i,j-1}+2U_{i+1,j-1}-20U_{ij}+8U_{i+1,j}+4U_{i,j+1}+2U_{i+1,j+1} \right) \notag	\\-\frac{1}{12}\left(V_{i,j-1}-V_{i-1,j}+8V_{ij}+3V_{i+1,j}+V_{i,j+1} \right)=\frac{2}{h}(g_N)_{j}. \label{eq-V-bound1}
\end{align}
 
However, on the left boundary $x=x_l$, the index $i=0$ and $v_{i-1,j}$ is out of domain $\Omega$.
 Meanwhile, $v$ is unknown, the previous treatment in \cite{li-pan-hoc21} is not practical. 
 Hence, we propose a more appropriate method to eliminate $v_{i-1,j}$, which does not destroy the compactness and accuracy.
 
 At first, we have Poisson equations $\Delta v=f$ hold at every grid point $(x_i,y_j)$. After using the second order accurate finite difference discretization, there is
 \begin{equation}
 \begin{aligned}
 	\frac{1}{h^2}\left(v(x_i,y_{j-1})+v(x_{i-1},y_j)-4v(x_i,y_j)+v(x_{i+1},y_j)+v(x_i,y_{j+1})\right)\\ = f(x_i,y_i)+O(h^2).
 \end{aligned}
 \end{equation}
 
 Here, we get an $O(h^4)$ approximations to $v(x_{i-1},y_j)$ as follows,
   \begin{equation}
  \begin{aligned}
 	v(x_{i-1},y_j) = h^2f(x_i,y_j)-\left(v(x_i,y_{j-1})-4v(x_i,y_j)+v(x_{i+1},y_j)+v(x_i,y_{j+1})\right) \\+O(h^4).
  \end{aligned}
   \end{equation}
 If $V_{ij}$ is also a high order approximation to $v(x_i,y_j)$,
 we also have
  \begin{equation}
  \begin{aligned}
 	V_{i-1,j} = h^2f_{ij}-\left(V_{i,j-1}-4V_{ij}+V_{i+1,j}+V_{i,j+1}\right) +O(h^\mu), \label{eq-Vappro}
  \end{aligned}
   \end{equation}
 where $\mu$ depends on ${V_{ij}}$'s accuracy and is less than or equal to 4.
 
When substituting \eqref{eq-Vappro} into \eqref{eq-V-bound1}, we can get the high order and compact discrete formulations
 \begin{align}
\frac{1}{6h^2}\left(4U_{i,j-1}+2U_{i+1,j-1}-20U_{ij}+8U_{i+1,j}+4U_{i,j+1}+2U_{i+1,j+1} \right) \notag	\\-\frac{1}{12}\left(2V_{i,j-1}+4V_{ij}+4V_{i+1,j}+2V_{i,j+1} \right)=\frac{2}{h}(g_N)_{j}-\frac{h^2}{12}f_{ij}. \label{eq-V-bound2}
\end{align}
\begin{remark}
	The truncation error of the finite difference formulation \eqref{eq-V-bound2} is $O(h^3)$. The analytical expression is $T_{ij} = (h^3/36u_{xyyyy}-7h^3/180u_{xxxxx})\big\vert_{(x_i,y_j)}$.	
\end{remark}
Now, we derive discrete boundary conditions for $v$'s Poisson problems at all boundary grid points.

 When all discrete equations for $v$ are combined, we can obtain 
 \begin{equation}
 	C{U}+D_h{V}=F_2.\label{eq-Dis-vPoiss}
 \end{equation}
 
 By solving the final linear system 
 \begin{equation}
 	\mathcal{L}_h\mathcal{U}=\begin{pmatrix}
 		A_h&B\\
 		C&D_h
 	\end{pmatrix}\begin{pmatrix}
 		\mathcal{U}\\
 		\mathcal{V}
 	\end{pmatrix} = \begin{pmatrix}
 		F_1\\
 		F_2
 	\end{pmatrix},\label{eq-4th-2d}
 \end{equation}
 the approximate values of $u$ and $v$ can be obtained. Furthermore, in numerical experiments, it has been confirmed that the conditioning number of the coefficient matrix $\mathcal{L}_h$ increases at a rate of $O(h^{-2})$. 

\begin{theorem}
	Let $\left\{U_{ij}, V_{ij}\vert 0\leq i\leq N_x,0\leq j\leq N_y\right\}$ and $\left\{u_{ij},v_{ij}\vert 0\leq i\leq N_x,0\right.$ $\left.\leq j\leq N_y\right\}$ be soultion of the finite difference scheme \eqref{eq-4th-2d} and the biharmonic problem \eqref{eq-vP}-\eqref{eq-vB2}, respectively. There is  constants $C_u$, $C_v$ independent of $h$ such that
	\begin{equation*}
		\vert\vert{U_{ij}-u_{ij}}\vert\vert_\infty \leq C_u h^4,\qquad \vert\vert{V_{ij}-v_{ij}}\vert\vert_{L_2} \leq C_v h^4.
	\end{equation*}
\end{theorem}

\section{Fourth order compact schemes for Biharmonic problem in three dimension}\label{sec:3d}
In this section, we focus on three-dimensional biharmonic problems. After introducing the intermediate variables of the Laplace of the original variables, we also obtain the two Poisson system
\begin{align}
\begin{cases}
	\Delta u =v,   \;\;\text{in}\; \Omega \\
	u = g_D,   \;\;\text{on}\; \partial\Omega
\end{cases},\qquad
\begin{cases}
	\Delta v =f,   \;\;\text{in}\; \Omega \\
	\partial_n u = g_N,   \;\;\text{on}\; \partial\Omega_1 \\
	\Delta u = g_L,   \;\;\text{on}\; \partial\Omega_2,
\end{cases}
\end{align}\label{eq-bi-3d}
where the cubic domain $\Omega$ is $\left[x_l,x_r \right]\times\left[y_l,y_r \right]\times\left[z_l,z_r \right]$. 
To derive the fourth order scheme clearly, we also assume that $\partial\Omega_1 = \left\{(x_l,y,z)\vert y_l<y<y_r,z_l<z<z_r\right\}$.

Similarly, we use a uniform grid $\Omega_h$ to cover the domain $\Omega$, which is defined by  
\begin{equation}
\begin{aligned}
	&x_i = x_l+ih,\;\;i=0,\cdots,N_x;\qquad y_j=y_l+jh,\;\;j=0,\cdots,N_y;\\
    &z_k=z_l+kh,\;\;j=0,\cdots,N_z, \label{Grid3D}
\end{aligned}
\end{equation}
and $h=(x_r-x_l)/N_x=(y_r-y_l)/N_y=(z_r-z_l)/N_z$.

For both $u$'s Poisson and $v$'s Poisson equation at an interior grid point $(x_i,y_j,z_k)$, we can apply the classical fourth order compact finite difference scheme. Then, we obtain the $u$'s discrete Poisson equations
\begin{equation}
\begin{aligned}
&-\frac{4}{h^2} U_{ijk} + \frac{1}{3h^2} \left( U_{i+1,j,k} + U_{i-1,j,k} + U_{i,j+1,k} + U_{i,j-1,k} + U_{i,j,k+1} + U_{i,j,k-1} \right) \\
&\quad + \frac{1}{6h^2} \left( U_{i+1,j+1,k} + U_{i+1,j-1,k} + U_{i-1,j+1,k} + U_{i-1,j-1,k} + U_{i+1,j,k+1} + U_{i+1,j,k-1} \right.\\ 
&\quad\left.+ U_{i-1,j,k+1}+ U_{i-1,j,k-1} + U_{i,j+1,k+1} + U_{i,j-1,k+1} + U_{i,j+1,k-1} + U_{i,j-1,k-1} \right) \\
&= \frac{1}{12} \left( 6 v_{ijk} + v_{i+1,j,k} + v_{i-1,j,k} + v_{i,j+1,k} + v_{i,j-1,k} + v_{i,j,k+1} + v_{i,j,k-1} \right), 
\end{aligned}\label{eq-u-int-3D}
\end{equation}
and $v$'s discrete Poisson equations
\begin{equation}
\begin{aligned}
&-\frac{4}{h^2} V_{ijk} + \frac{1}{3h^2} \left( V_{i+1,j,k} + V_{i-1,j,k} + V_{i,j+1,k} + V_{i,j-1,k} + V_{i,j,k+1} + V_{i,j,k-1} \right) \\
&\quad + \frac{1}{6h^2} \left( V_{i+1,j+1,k} + V_{i+1,j-1,k} + V_{i-1,j+1,k} + V_{i-1,j-1,k} + V_{i+1,j,k+1} + V_{i+1,j,k-1}\right.\\ 
&\quad\left. + V_{i-1,j,k+1}+ U_{i-1,j,k-1} + U_{i,j+1,k+1} + U_{i,j-1,k+1} + U_{i,j+1,k-1} + U_{i,j-1,k-1} \right) \\
&= \frac{1}{12} \left( 6 f_{ijk} + f_{i+1,j,k} + f_{i-1,j,k} + f_{i,j+1,k} + f_{i,j-1,k} + f_{i,j,k+1} + f_{i,j,k-1} \right),
\end{aligned}\label{eq-v-int-3D}
\end{equation}
where $U_{ijk}$ and $V_{ijk}$ are the finite difference solutions to biharmonic problems, and $f_{ijk}$ is the value of the right term at the grid points.

Since the values of $u$ on the boundary are known, we get the $u$'s discrete Poisson problems in matrix form
\begin{equation}
	A_h\mathcal{U}+B\mathcal{V}=F_1, \label{eq-Dis-uPoiss-3D}
\end{equation}

Compared to two dimensional cases, three kinds of boundary grid points need to be considered: corner grid points (e.g. $(x_l,y_l,z_l)$), edge grid points (e.g. $(x_l,y_l,z_k)$), face grid points (e.g. $(x_l,y_j,z_k)$).

For the first two situations, using the Dirichlet boundary conditions of $u$, the values of $v$ can be computed directly, which avoids constructing the finite difference discretization.
Such as corner grid points $(x_l,y_l,z_l)$, we can differentiate the $u=g_D$ from three different directions (parallel to the $x$-axis, $y$-axis, $z$-axis) to obtain $u_{xx}(x,y_l,z_l),u_{yy}(x_l,y,z_l),u_{zz}(x_l,y_l,z)$. 
Further, the value of $v=\Delta u$ at corner gird point $(x_l,y_l,z_l)$ can be computed. 
The edge $(x_l,y_l,z)$ is intersection of the face $(x_l,y,z)$  and $(x,y_l,z)$. 
We want to get $v$ on edge $(x_l,y_l,z)$ and know $u(x_l,y,z)$ and $u(x,y_l,z)$ from the given conditions.
Similarly, by differentiating the boundary conditions of $u$, we can derive $v$ on the edge.

\begin{remark}
  If $\partial\Omega_1=\partial\Omega$ and $\partial\Omega_2=\emptyset$, we can use the same derivation to get $v$ on all corner grid points and edge grid points. There do not exist any differences. Through these treatments, it is equivalent to getting the $v$'s Dirichlet boundary conditions on the corner and edge.
\end{remark}

The main challenge is constructing the finite difference discretization on the face grid points $(x_l,y_j,z_k)$.

Different from two dimensional cases, we detail a modified methodology of \cite{li-pan-hoc21} to get a discrete formulation for $v$ on the boundary in three dimensions.
The basis is that we need to treat $v$ as the source term of $u$'s Poisson equation and combine the Neumann boundary condition $\partial_n u=g_N$ and the load $f$.

Assume that a compact finite difference scheme can be written as
\begin{equation}
	\begin{aligned}
		\sum_{i_e=0}^1\sum_{j_e=-1}^1\sum_{k_e=-1}^1 \alpha_{i_e,j_e,k_e} U_{i+i_e,j+j_e,k+k_e}-\sum_{i_e=0}^1\sum_{j_e=-1}^1\sum_{k_e=-1}^1 \beta_{i_e,j_e,k_e} V_{i+i_e,j+j_e,k+k_e}\\= 
	 f(x_{i},y_{j},z_{k})
	+ \sum_{j_e=-1}^1\sum_{k_e=-1}^1 \gamma_{j_e,k_e} g_N(y_{j+j_e},z_{k+k_e}) \label{eq-bound-3d-fd}
	\end{aligned}
\end{equation}
where $i=0$ throughout this section, $\alpha_{i_e,j_e,k_e}, \beta_{i_e,j_e,k_e}(i_e=0,1,j_e=-1,0,1,k_e=-1,0,1)$, and $\gamma_{j_e,k_e} (j_e=-1,0,1,k_e=-1,0,1)$ are undermined coefficients.
We leave the index $i$ in the formulas so that the derivation can be applied to other boundaries when $\partial\Omega_1 =\partial\Omega$. 
Thus, we have 45 coefficients, including 18 coefficients for $ U_{i+i_e,j+j_e,k+k_e}$, 18 coefficients for $V_{i+i_e,j+j_e,k+k_e}$, 9 coefficients for $g_N(y_{j+j_e},z_{k+k_e})$. 

Denote the truncation error of the scheme \eqref{eq-bound-3d-fd} at the boundary grid points $(x_i,y_j,z_k)$ as follows,
\begin{equation}
	\begin{aligned}
	&T_{ijk} = 	\sum_{i_e=0}^1\sum_{j_e=-1}^1\sum_{k_e=-1}^1 \alpha_{i_e,j_e,k_e} u(x_{i+i_e},y_{j+j_e},z_{k+k_e})- f(x_{i},y_{j},z_{k})\\
	&-\sum_{i_e=0}^1\sum_{j_e=-1}^1\sum_{k_e=-1}^1 \beta_{i_e,j_e,k_e} v(x_{i+i_e},y_{j+j_e},z_{k+k_e}) 
	-\sum_{j_e=-1}^1\sum_{k_e=-1}^1 \gamma_{j_e,k_e} g_N(y_{j+j_e},z_{k+k_e}).\label{eq-trunc-3d}
	\end{aligned}
\end{equation}
We want to determine the coefficients in \eqref{eq-bound-3d-fd} such that the local truncation errors \eqref{eq-trunc-3d} are zero or $O(h^4)$ for any fourth order polynomials. 
To derive a linear system for the solution of the coefficients, we will expand the term $u(x_{i+i_e},y_{j+j_e},z_{k+k_e}), v(x_{i+i_e},$ $y_{j+j_e},z_{k+k_e})$ at the grid point $(x_i,y_j,z_k)$ and $ g_N(y_{j+j_e},z_{k+k_e})$ at $(y_j,z_k)$.

For the term $u(x_{i+i_e},y_{j+j_e},z_{k+k_e})$, apply the multivariable Taylor expansion at grid points $(x_i,y_j,z_k)$ up to all fourth order partial derivatives, 
\begin{equation}
\small	u(x_{i+i_e},y_{j+j_e},z_{k+k_e}) = \sum_{0\leq e_i+e_j+e_k\leq4}\frac{1}{e_i!e_j!e_k!}\frac{\partial^{e_i+e_j+e_k}u}{\partial x^{e_i}\partial y^{e_j}\partial z^{e_k}}\Bigg\vert_{x_i,y_j,z_k}h_{i_e}^{e_i}h_{j_e}^{e_j}h_{k_e}^{e_k} +O(h^5), \label{eq-u-expan}
\end{equation}
where $h_{i_e}=i_eh, h_{j_e}=j_eh$ and $h_{k_e}=k_eh$.

Apply the Taylor expansion to the term $v(x_{i+i_e},y_{j+j_e},z_{k+k_e})$ up to all second order partial derivatives can get
\begin{equation}
\small		v(x_{i+i_e},y_{j+j_e},z_{k+k_e}) = \sum_{0\leq e_i+e_j+e_k\leq2}\frac{1}{e_i!e_j!e_k!}\frac{\partial^{e_i+e_j+e_k}v}{\partial x^{e_i}\partial y^{e_j}\partial z^{e_k}}\Bigg\vert_{(x_i,y_j,z_k)}h_{i_e}^{e_i}h_{j_e}^{e_j}h_{k_e}^{e_k} +O(h^3).\label{eq-v-expan}
\end{equation}
And by utilizing $v=\Delta u$ and differentiating $v$, we can also obtain the derivatives of $v$ expressed by $u$,
\begin{equation}
	\frac{\partial^{e_i+e_j+e_k}v}{\partial x^{e_i}\partial y^{e_j}\partial z^{e_k}} = \frac{\partial^{e_i+e_j+e_k+2}u}{\partial x^{e_i+2}\partial y^{e_j}\partial z^{e_k}}+\frac{\partial^{e_i+e_j+e_k+2}u}{\partial x^{e_i}\partial y^{e_j+2}\partial z^{e_k}}+\frac{\partial^{e_i+e_j+e_k+2}u}{\partial x^{e_i}\partial y^{e_j}\partial z^{e_k+2}}
\end{equation}
where $0\leq e_i+e_j+e_k\leq2$, such as, $v_{x} = u_{xxx}+u_{xyy}+u_{xzz}$.

For the term $g_N(y_{j+j_e},z_{z+z_k})$, we apply the Taylor expansion up to third order derivatives
\begin{equation}
\small		g(y_{j+j_e},z_{k+k_e}) = \sum_{0\leq e_j+e_k\leq2}\frac{1}{e_j!e_k!}\frac{\partial^{e_j+e_k}g_N}{\partial y^{e_j}\partial z^{e_k}}\Bigg\vert_{(y_j,z_k)}h_{j_e}^{e_j}h_{k_e}^{e_k} +O(h^4). \label{eq-NC-expan}
\end{equation}
Note that 
\begin{equation}
	\partial_n u\Big\vert_{x=x_l} = -u_x\Big\vert_{x=x_l} =g_N(y,z).
\end{equation}
Hence, in the expansion \eqref{eq-NC-expan}, all partial derivatives of $g_N$ can be substituted by $u$ and its partial derivatives,
\begin{equation}
	\frac{\partial^{e_j+e_k}g_N}{\partial y^{e_j}\partial z^{e_k}} = -\frac{\partial^{e_j+e_k+1}u}{\partial x\partial y^{e_j}\partial z^{e_k}}, \qquad 0\leq e_j+e_k\leq2.
\end{equation}

After the expansions of all involved terms and adding $f= \Delta^2 u$ at grid points $(x_i,y_j,z_k)$, we can rewrite the local truncation error as
\begin{equation}
	T_{ijk} = \sum_{0\leq e_i+e_j+e_k\leq 4}L_{e_i,e_j,e_k}\frac{\partial^{e_i+e_j+e_k}u}{\partial x^{e_i}\partial y^{e_j}\partial z^{e_k}}+O(\vert\vert{\boldsymbol \alpha}\vert\vert_\infty h^5+\vert\vert{\boldsymbol \beta}\vert\vert_\infty h^3+\vert\vert{\boldsymbol \gamma}\vert\vert_\infty h^4),
\end{equation}
where $L_{e_i,e_j,e_k}$ are the results after we collect terms \eqref{eq-u-expan},\eqref{eq-v-expan} and \eqref{eq-NC-expan}, and are the linear combination of the coefficients $\alpha_{i_e,j_e,k_e}, \beta_{i_e,j_e,k_e},\gamma_{j_e,k_e}$. 
Note that $\vert\vert{\boldsymbol \alpha}\vert\vert_\infty = \max_{0\leq i_e\leq 1,-1\leq j_e\leq 1,-1 \leq k_3 \leq 1}\{\vert\alpha_{i_e,j_e,k_e}\vert\}$, and so on.
Here, we want the truncation error $T_{ijk}$ to be zeros or $O(h^4)$ for all fourth order polynomials, or $x^{e_i}y^{e_j}z^{e_k}$, $0\leq e_i+e_j+e_k\leq 4$. 
In this way, by matching the terms of the coefficients of $u,u_x,u_y,u_z,\cdots,$ $u_{xxxx},u_{xxxy},\cdots,u_{yyyy},u_{yyyz},\cdots,u_{zzzz}$, we can derive 35 linear equations whose solution are the undetermined coefficients $\alpha_{i_e,j_e,k_e}, \beta_{i_e,j_e,k_e}, \gamma_{j_e,k_e}$.
\begin{equation}
	\begin{aligned}
	&\sum_{i_e=0}^1\sum_{j_e=-1}^1\sum_{k_e=-1}^1\alpha_{i_e,j_e,k_e} \frac{h_{i_e}^{e_i}h_{j_e}^{e_j}h_{k_e}^{e_k}}{e_i!e_j!e_k!}-\sum_{i_e=0}^1\sum_{j_e=-1}^1\sum_{k_e=-1}^1\beta_{i_e,j_e,k_e}\\
	&\left(  \frac{h_{i_e}^{e_i-2}h_{j_e}^{e_j}h_{k_e}^{e_k}}{(e_i-2)!e_j!e_k!}H(e_i-1) +  \frac{h_{i_e}^{e_i}h_{j_e}^{e_j-2}h_{k_e}^{e_k}}{e_i!(e_j-2)!e_k!} H(e_j-1) +  \frac{h_{i_e}^{e_i}h_{j_e}^{e_j}h_{k_e}^{e_k-2}}{e_i!e_j!(e_k-2)!} H(e_k-1) \right)\\
	&-\sum_{j_e=-1}^1\sum_{k_e=-1}^1 \gamma_{j_e,k_e} \frac{h_{j_e}^{e_j}h_{k_e}^{e_k}}{e_j!e_k!}\chi(e_i-1)= \chi(e_i-4) +  \chi(e_j-4) + \chi(e_k-4) \\
	&+ 2\left(\chi(e_i-2)\chi(e_j-2)\chi(e_k) + \chi(e_i-2)\chi(e_j)\chi(e_k-2) + \chi(e_i)\chi(e_j-2)\chi(e_k-2)\right),
	\end{aligned}
\end{equation}
where $\chi(x)$ is Characteristic function and $H(x)$ is Heaviside function defined as 
\begin{equation}
	\chi(x)=\begin{cases}
		1, \quad x=0,\\
		0, \quad x\neq0,
	\end{cases}\qquad
	H(x)=\begin{cases}
		1, \quad x\geq0,\\
		0, \quad x\leq0.
	\end{cases}
\end{equation}

When the parameters $e_i, e_j, e_k$ set by a range of $0$ to $4$ and satisfy $0\leq e_i+e_j+e_k\leq4$,  we can obtain 35 linear equations containing 45 unknowns. This is an underdetermined system of equations. In general there is an infinite number of solutions. Using the Symbolic package of Matlab, we have found a set of coefficients analytically for the $v$'s discretization on the boundary $x=x_l$. The specific discrete formulations are shown as follows,
\begin{equation}
	\begin{aligned}
		\frac{1}{6h^2}\left(U_{i,j-1,k-1}+2U_{i,j,k-1}+U_{i,j+1,k-1}+2U_{i,j-1,k} -24U_{ijk}+2U_{i,j+1,k} \right.\\
	\left.+U_{i,j-1,k+1}+2U_{i,j,k+1}+U_{i,j+1,k+1}\right)+\frac{1}{6h^2}\left(2U_{i+1,j-1,k}+U_{i+1,j,k-1}+ 4U_{i+1,j,k}\right.\\
	\left.+U_{i+1,j+1,k}+U_{i+1,j,k+1}\right)-\frac{1}{12}\left( V_{i,j-1,k}+V_{i,j,k-1}+ 4V_{ijk}+V_{i,j+1,k}+V_{i,j,k+1}\right.\\
	\left. + V_{i+1,j-1,k}+V_{i+1,j,k-1}+V_{i+1,j+1,k}+V_{i+1,j,k+1} \right)=-\frac{h^2}{12}f(x_i,y_j,z_k)-\frac{2}{h}g_N(y_j,z_k).
	\end{aligned}\label{eq-v-bound-3d}
\end{equation}
Here, we divide $(-12/h^2)$ in the above discrete formula such that the coefficient of the load $f$ does not equal to 1 at the grid point $(x_i,y_j,z_k)$.
When we combine the discrete formulas  \eqref{eq-v-int-3D}, \eqref{eq-v-bound-3d} and the Dirichlet boundary conditions of $v$, the matrix form of $v$'s Poisson problem can be wrriten as
 \begin{equation}
 	C{U}+D_h{V}=F_2.\label{eq-Dis-vPoiss-3d}
 \end{equation}
 
 The final discrete coupled system of the biharmonic \eqref{eq-bi-3d} can be obtained by
 \begin{equation}
 	\mathcal{L}_h\mathcal{U}=\begin{pmatrix}
 		A_h&B\\
 		C&D_h
 	\end{pmatrix}\begin{pmatrix}
 		\mathcal{U}\\
 		\mathcal{V}
 	\end{pmatrix} = \begin{pmatrix}
 		F_1\\
 		F_2
 	\end{pmatrix}.\label{eq-4th-3d}
 \end{equation}
Note that in 3D, similar to 2D, the unknowns only contain $u$ and $v=\Delta u$. This does not increase with dimension compared to the very popular Stephenson scheme.
 Next, we discuss the convergence of the fourth order compact scheme for the biharmonic problem as summarized in the following theorem.
\begin{theorem}
	Let $\left\{U_{ijk}, V_{ijk}\vert 0\leq i\leq N_x,0\leq j\leq N_y,0\leq k\leq N_z\right\}$ and $\left\{u_{ijk},\right.$$\left.v_{ijk}\vert 0\leq i\leq N_x,0\leq j\leq N_y,0\leq k\leq N_z\right\}$ be soultion of the finite difference scheme \eqref{eq-4th-3d} and the biharmonic problem \eqref{eq-bi-3d}, respectively. There is  constants $C_u$, $C_v$ independent of $h$ such that
	\begin{equation*}
		\vert\vert{U_{ij}-u_{ij}}\vert\vert_\infty \leq C_u h^4,\qquad \vert\vert{V_{ij}-v_{ij}}\vert\vert_{L_2} \leq C_v h^4.
	\end{equation*}
\end{theorem}

\section{Numerical experiments} \label{sec:NE}
We have tested the proposed fourth order compact scheme for polynomials $P^n(x, y)$ in 2D, $P^n(x, y,z)$ in 3D, $n \leq 4$; the primary variables $\mathcal{U}$ can be accurate to $\varepsilon$cond($A_h$),  where $\epsilon \sim 1e-16$ is the machine precision and $A_h$ is the coefficient matrix of the finite difference equations of a Poisson equation not the discrete coupled system of biharmonic problems $\mathcal{L}_h$. 
However, the intermediate variables $\mathcal{V}$ are only accurate to $\varepsilon$cond($\mathcal{L}_h$). 
Next, we test some constructed examples with genuine nonlinear solutions, one of which is relatively smooth and the other can be oscillatory.

\subsection{Two dimensional examples}
In this subsection, we test the fourth order scheme for the biharmonic with the Dirichlet boundary conditions in 2D. And the constructed examples come from published papers. 

\begin{example} \label{ex-s-2d}
	The following is an example with a smooth solution:
	\begin{align*}
		u(x,y) &=x^2+y^2-xe^x\cos(y), \qquad (x,y)\in (0,1)^2,\\
		\partial_n u(x,y) &= \left[2x - e^x\cos(y) - xe^x\cos(y), 2y + xe^x\sin(y)\right]\cdot \vec{n}, \qquad (x,y) \in \partial\Omega_1,\\
		\Delta u(x,y) &=  4 - 2e^x\cos(y),\qquad (x,y) \in \partial\Omega_2,
	\end{align*}
where $\vec{n}$ are out normal unit vectors,  $\partial \Omega_1\cup\partial\Omega_2=\partial\Omega$ and the load $f$ can be calculated analytically. 
\end{example}
Two kinds of boundary conditions are tested containing: $\partial \Omega_1=\partial\Omega$, Dirichlet boundary conditions of the first kind; $\partial\Omega_2=\{(0,y)\vert 0\leq y\leq 1\}$, Dirichlet boundary conditions of the mix.

In table \ref{tab_2d_ex1}, we list the error in infinity norm and convergent order when two different kinds of BCs are prescribed. 
We can see a clear fourth order for Dirichlet BCs of first kind. 
However, due to the dominant effect of round-off errors on the finest mesh, the convergent order becomes super-third for Dirichlet boundary conditions of mixed type.

\begin{example}\label{ex-o-2d}
	The following is an example with an oscillatory solution:
	\begin{align*}
		u(x,y) &=\sin(k_1 x)\cos(k_2 y);, \qquad (x,y)\in (0,1)^2,\\
		\partial_n u(x,y) &= \left[ k_1\cos(k_1 x)\cos(k_2 y), -k_2\sin(k_1 x)\sin(k_2 y) \right]\cdot \vec{n}, \qquad (x,y) \in \partial\Omega_1,\\
		\Delta u(x,y) &=-(k_1^2+ k_2^2)\sin(k_1 x)\cos(k_2 y) ,\qquad (x,y) \in \partial\Omega_2,
	\end{align*}
where $\vec{n}$ are out normal unit vectors,  $\partial \Omega_1\cup\partial\Omega_2=\partial\Omega$ and the load $f$ can be also calculated analytically. 
\end{example}
In this example, we can choose $k_1$ and $k_2$ to make the solution more oscillatory.

Table \ref{tab_2d_ex2_1} and Table \ref{tab_2d_ex2_2} show the grid refinement analysis for Example \ref{ex-s-2d} with two pairs of $k_1$ and $k_2$. In Table \ref{tab_2d_ex2_1}, we chose $k_1=25$ and $k_2=5$. Regardless of whether Dirichlet boundary conditions of the first kind or mixed type are prescribed, the fourth-order convergence is verified.
In Table \ref{tab_2d_ex2_2}, we set $k_1=5,k_2=50$ such that the solution is more oscillatory. We observe that the numerical solution is highly accurate on a $1024\times1024$ mesh, but it is challenging for low order schemes. And convergent order attains fourth order.

\begin{figure}[htbp]
\centering
\includegraphics[width=0.7\textwidth]{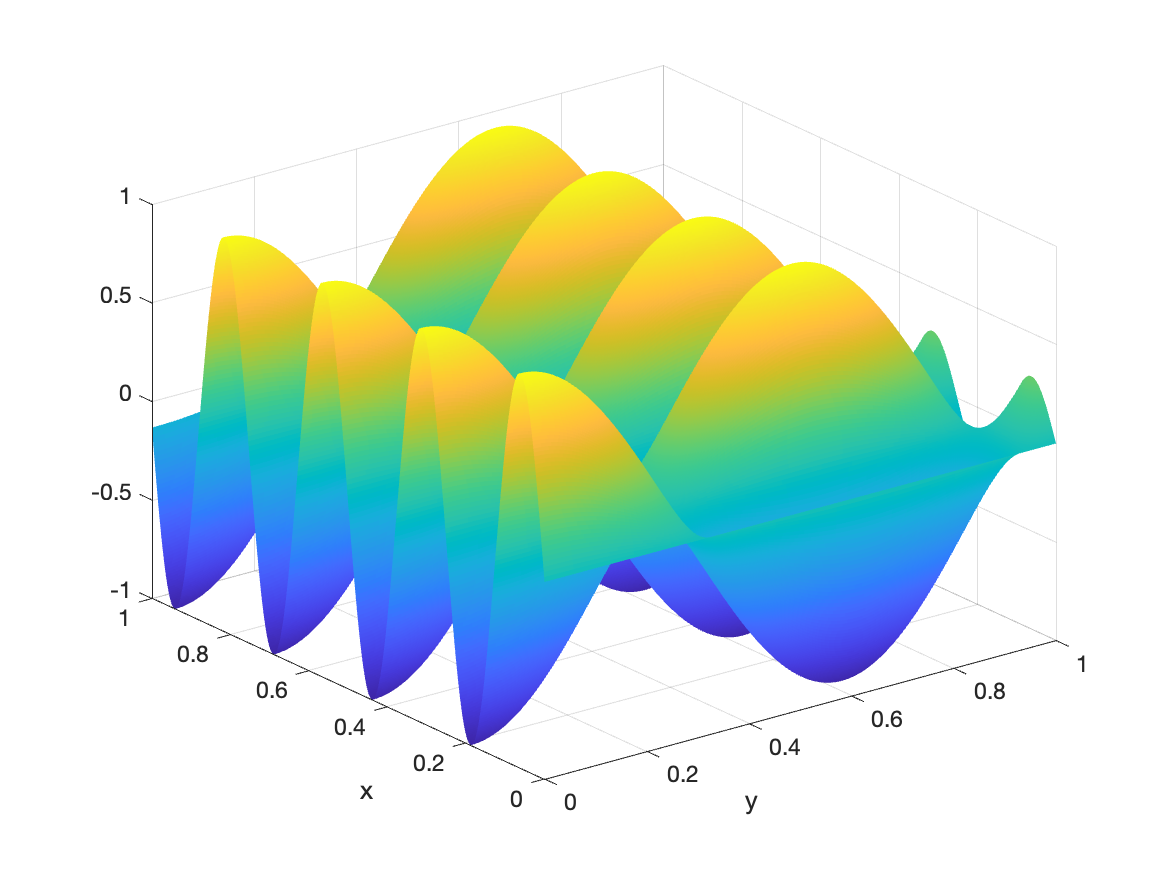}
\caption{The numerical solutions plot for Exampel \ref{ex-o-2d} with $k_1=25,k_2=5$ on $1024\times$1024 mesh. } \label{fig_25_5}
\end{figure}

\begin{figure}[htbp]
\centering
\includegraphics[width=0.7\textwidth]{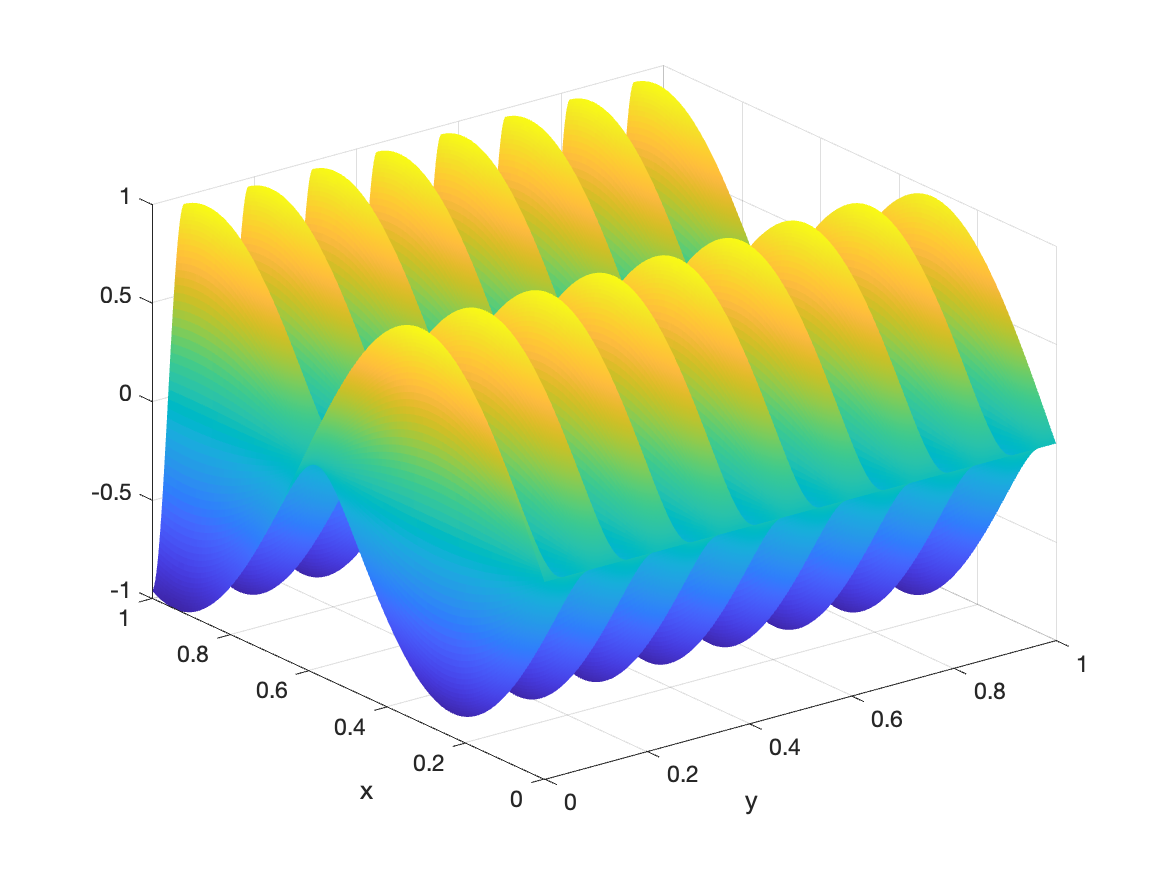}
\caption{The numerical solutions plot for Exampel \ref{ex-o-2d} with $k_1=5,k_2=50$ on $1024\times$1024 mesh. } \label{fig_5_50}
\end{figure}

In Table \ref{tab_2d_cond}, we show the conditioning number of the coefficient matrix of the proposed scheme and second order 13-point scheme.
Despite the introduction of an extra unknown in the proposed scheme, its conditioning number is lower than that of the 13-point scheme.
Furthermore, the increasing rate is approximately $O(h^{-2})$.
Figure \ref{fig_25_5} and \ref{fig_5_50} show the numerical solution plot of Example \ref{ex-o-2d}.

\begin{table}[htbp]\label{tab_2d_ex1}
\begin{center}
\caption{Grid refinement analysis of the fourth order compact scheme for Example \ref{ex-s-2d}. }
\centering
\begin{tabular}{ccccc}
  \hline
  \multirow{2}{*}{$N$} & \multicolumn{2}{c}{Dirichelt BCs of first kind} & \multicolumn{2}{c}{Dirichlet BCs of mix}  \\
  & $\| E \|_\infty$  & Order    & $\|E\|_{\infty}$  & Order  \\
  \hline
64 & 8.53e-05 & -- & 1.72e-05 & --  \\
128 & 5.37e-06 & 3.99 & 1.08e-06 & 4.00 \\
256 & 3.37e-07 & 4.00 & 6.73e-08 & 4.00\\
512 & 2.11e-08 & 4.00 & 4.24e-09 & 3.99\\
1024 &  1.31e-09 & 4.01 & 4.07e-10 & 3.38\\
\hline
\end{tabular}
\end{center}
\end{table}

\begin{table}[htbp]\label{tab_2d_ex2_1}
\begin{center}
\caption{Grid refinement analysis of the fourth order compact scheme for Example \ref{ex-o-2d} with $k_1=25,k_2=5$. }
\centering
\begin{tabular}{ccccc}
  \hline
  \multirow{2}{*}{$N$} & \multicolumn{2}{c}{Dirichelt BCs of first kind} & \multicolumn{2}{c}{Dirichlet BCs of mix}  \\
  & $\| E \|_\infty$  & order    & $\|E\|_{\infty}$  & order  \\
  \hline
64 &5.12e-04& -- & 7.21e-04 & -- \\
128 &3.23e-05& 3.99 & 4.53e-05 & 3.99\\
256 &2.02e-06& 4.00 & 2.83e-06 & 4.00\\
512 & 1.26e-07&  4.00 & 1.77e-07 & 4.00 \\
1024 & 7.90e-09&4.00 & 1.11e-08 & 4.00\\
\hline
\end{tabular}
\end{center}
\end{table}

\begin{table}[htbp]\label{tab_2d_ex2_2}
\begin{center}
\caption{Grid refinement analysis of the fourth order compact scheme for Example \ref{ex-o-2d} with $k_1=5,k_2=50$. }
\centering
\begin{tabular}{ccccc}
  \hline
  \multirow{2}{*}{$N$} & \multicolumn{2}{c}{Dirichelt BCs of first kind} & \multicolumn{2}{c}{Dirichlet BCs of mix}  \\
  & $\| E \|_\infty$  & order    & $\|E\|_{\infty}$  & order  \\
  \hline
64 & 1.16e-02 & -- & 3.05e-02 & -- \\
128 & 6.38e-04 & 4.19 & 1.87e-03 & 4.03\\
256 & 3.68e-05 & 4.12 & 1.16e-04 & 4.01\\
512 & 2.20e-06 & 4.06 & 7.27e-06 & 4.00\\
1024 & 1.36e-07 & 4.02 & 4.54e-07 & 4.00  \\
\hline
\end{tabular}
\end{center}
\end{table}

\begin{table}[htbp]\label{tab_2d_cond}
\begin{center}
\caption{A comparison in conditioning number of the coefficient matrix in two dimension. }
\centering
\begin{tabular}{ccccc}
  \hline
  \multirow{2}{*}{$N$} & \multicolumn{4}{c}{Dirichelt BCs of 1st kind}   \\
  & Our scheme  & Rate     & 13-point scheme & Rate \\
  \hline
128 & 1.62e+07 & -- & 2.17e+07 & --\\
256 & 6.50e+07 & 4.00 & 3.48e+08 & 15.99\\
512 & 2.60e+08 & 4.00 & 5.57e+09  & 16.00\\
1024 & 1.04e+09 & 4.00 & 8.90e+10 & 16.00\\
2048 &  4.16e+09 & 4.00 & 1.42e+12 & 16.00\\
\hline
\end{tabular}
\end{center}
\end{table}

\begin{example}\label{ex-stokes}
	Application to incompressible Stokes flow on $\Omega=[0,1]^2$ in stream-vorticity from. To verify the fourth order convergence, we modify the Lid-driven cavity flow. After modification, the boundary conditions are prescribed by
	\begin{equation}
		\begin{aligned}
		&u=0,\qquad (x,y)\in \partial\Omega, \\
		&[-u_y,u_x]=\begin{cases}
			[0,0], x=0,1;y=0,\\
		    [x^6(x-1)^6,0], y=1,	
		\end{cases}	
		\end{aligned}
	\end{equation}
The source term $f$ is zero for the whole domain.
Here, $u$ represents the stream function in Stokes flow.
\end{example}

Since we cannot obtain analytical solutions for Example \ref{ex-stokes}, we consider the numerical solutions on a $1024\times 1024$ mesh as reference solutions. 
Based on the reference solutions, we get the grid refinement analysis in Table \ref{tab-stokes}. The convergence order is almost fourth order. 
We believe our scheme can be applied to Navier-Stokes problems in stream-vorticity form after some minor modifications.

\begin{table}[htbp]\label{tab-stokes}
\begin{center}
\caption{Grid refinement analysis of the fourth order compact scheme for Example \ref{ex-stokes}. }
\centering
\begin{tabular}{ccc}
  \hline
  $N$   & $\|E\|_{\infty}$  & Order  \\
  \hline
16 & 3.88e-07 & --   \\
32 & 3.81e-08 & 3.35   \\
64 & 2.27e-09 & 4.07\\
128 & 1.66e-10 & 3.78\\
256 &  1.17e-11 & 3.83 \\
\hline
\end{tabular}
\end{center}
\end{table}

\subsection{Three dimensional examples}

\begin{example}\label{ex-3d-1}
The following is an example with a smooth solution:
	\begin{align*}
		u(x,y) &=xyz\log(x+y+z+1);, \qquad (x,y)\in (0,1)^3,\\
		\partial_n u(x,y) &= \frac{xyz}{x + y + z + 1}+\left[yzlog(x + y + z + 1),\right.\\
		&\left. xzlog(x + y + z + 1), xylog(x + y + z + 1) \right]\cdot \vec{n}, \qquad (x,y) \in \partial\Omega_1,\\
		\Delta u(x,y) &=\frac{2x^2y + 2x^2z + 2xy^2 + 3xyz + 2xy + 2xz^2 + 2xz + 2y^2z + 2yz^2 + 2yz}{(x + y + z + 1)^2}
 ,\qquad (x,y) \in \partial\Omega_2,
	\end{align*}
	where $\vec{n}$ are out normal unit vectors,  $\partial \Omega_1\cup\partial\Omega_2=\partial\Omega$ and the load $f$ can be calculated analytically. 
\end{example}

\begin{example}\label{ex-3d-2}
The following is an example with an oscillatory solution:
	\begin{align*}
		u(x,y) &=\sin(k_1 x)\cos(k_2 y);, \qquad (x,y)\in (0,1)^3,\\
		\partial_n u(x,y) &= \left[k_1\cos(k_1 x)\cos(k_2 y)\sin(k_3 z), -k_2\sin(k_1 x)\sin(k_2 y)\sin(k_3 z),\right.\\ 
		&\left. k_3\cos(k_2 y)\cos(k_3 z)\sin(k_1 x) \right]\cdot \vec{n}, \qquad (x,y) \in \partial\Omega_1,\\
		\Delta u(x,y) &=- (k_1^2+k_2^2+k_3^2 )\sin(k_1 x)\cos(k_2 y)\sin(k_3 z) ,\qquad (x,y) \in \partial\Omega_2,
	\end{align*}
	where $\vec{n}$ are out normal unit vectors,  $\partial \Omega_1\cup\partial\Omega_2=\partial\Omega$ and the load $f$ can be calculated analytically. 
\end{example}

For three three-dimensional problem, two kinds of boundary conditions are tested containing: $\partial \Omega_1=\partial\Omega$, Dirichlet boundary conditions of the first kind; $\partial\Omega_2=\{(0,y,z)\vert 0\leq y\leq 1,\vert 0\leq z\leq 1\}$, Dirichlet boundary conditions of the mix.

In Table \ref{tab_3d_s}, we simulated an example with smooth solutions. We observed that on a relatively fine mesh, the error almost reaches machine precision.

By controlling the parameters $k_1,k_2,k_3$, we also simulate an oscillatory solution whose results are listed in Table \ref{tab_3d_o}. The convergence order almost reaches fourth-order accuracy with the refined mesh.

Table \ref{tab_3d_cond} shows the conditioning number of our schemes and 27-points scheme \cite{FD3d27}. Our proposed scheme maintains a $O(h^{-2})$ rate, which is equal to two-dimensional cases. 

\begin{table}[htbp]\label{tab_3d_s}
\begin{center}
\caption{Grid refinement analysis of the fourth order compact scheme for Example \ref{ex-3d-1}. }
\centering
\begin{tabular}{ccccc}
  \hline
  \multirow{2}{*}{$N$} & \multicolumn{2}{c}{Dirichelt BCs of first kind} & \multicolumn{2}{c}{Dirichlet BCs of mix}  \\
  & $\| E \|_\infty$  & Order    & $\|E\|_{\infty}$  & Order  \\
  \hline
16 & 7.36e-08 & -- & 6.28e-08 & --\\
32 & 4.83e-09 & 3.93 & 4.25e-09 & 3.89 \\
64 & 3.06e-10 & 3.98 & 2.74e-10 & 3.96\\
128 & 1.93e-11 & 3.99 & 1.73e-11 & 3.98\\
256 & 1.21e-12 & 4.00 & 1.08e-11 & 4.00\\
\hline
\end{tabular}
\end{center}
\end{table}

\begin{table}[htbp]\label{tab_3d_o}
\begin{center}
\caption{Grid refinement analysis of the fourth order compact scheme for Example \ref{ex-3d-2} with $k_1=25,k_2=5,k_3=25$. }
\centering
\begin{tabular}{ccccc}
  \hline
  \multirow{2}{*}{$N$} & \multicolumn{2}{c}{Dirichelt BCs of first kind} & \multicolumn{2}{c}{Dirichlet BCs of mix}  \\
  & $\| E \|_\infty$  & Order    & $\|E\|_{\infty}$  & Order  \\
  \hline
16 & 1.36e-01 & -- & 1.19e-01 & --\\
32 & 1.21e-02 & 3.50 & 1.15e-02 & 3.37\\
64 & 9.31e-04 & 3.70 & 9.02e-04 & 3.67\\
128 & 5.98e-05 & 3.96&5.84e-05 & 3.95\\
256 & 3.77e-06 & 3.98 & 3.65e-06 & 4.00\\
\hline
\end{tabular}
\end{center}
\end{table}

\begin{table}[htbp]\label{tab_3d_cond}
\begin{center}
\caption{A comparison in conditioning number of the coefficient matrix in three dimensions. }
\centering
\begin{tabular}{ccccc}
  \hline
  \multirow{2}{*}{$N$} & \multicolumn{4}{c}{Dirichelt BCs of 1st kind}   \\
  & Our scheme  & Rate     & 13-point scheme & Rate \\
  \hline
16 &1.18e+06 & -- & 4.57e+05 & --\\
32 & 4.91e+06 & 4.16 & 6.68e+06 & 15.81\\
64 & 1.96e+07 & 3.99 & 1.09e+08 & 16.00\\
128 & 1.37e+08 & 4.00 & 1.74e+09 & 16.00\\
256 & 5.48e+08 & 4.00 & 2.74e+10 & 16.00\\
\hline
\end{tabular}
\end{center}
\end{table}

\section{Conclusion}
In this paper, we present a fourth-order compact finite difference scheme for solving biharmonic problems based on a coupled approach. 
This scheme has the advantage of not increasing the number of unknowns when moving from two to three dimensions, which addresses a drawback of Stephenson schemes. 
Additionally, the conditioning number of the discrete coefficient matrix increases at a rate of $O(h^{-2})$, providing flexibility in choosing a solver. 
We will also explore applying the scheme to nonlinear biharmonic problems and incompressible Navier-Stokes problems in stream-vorticity in future work.

\section*{Acknowledgements}

 Kejia Pan was supported by the National Natural Science Foundation of China (No. 42274101).
 Jin Li was supported by the Fundamental Research Funds for the Central Universities of Central South University  (No. 2022ZZTS0610).
Zhilin Li was partially supported by a Simons grant 633724.

%\section*{References}
\bibliographystyle{unsrt}

%\bibliography{references}
%\begin{thebibliography}{10}
%
%
%\end{thebibliography}

\end{document}